\documentclass[reqno]{amsart}

\usepackage{amsmath, amssymb, amsthm, epsfig}
\usepackage{hyperref, latexsym}
\usepackage{url}
\usepackage[mathscr]{euscript}
\usepackage{color}
\usepackage{harpoon}
\usepackage{url}

\usepackage{setspace}
\usepackage{refcheck} 
\norefnames
\nocitenames

\def\today{\ifcase\month\or
  January\or February\or March\or April\or May\or June\or
  July\or August\or September\or October\or November\or December\fi
  \space\number\day, \number\year}

\newtheorem{theorem}{Theorem}

\theoremstyle{definition}

\theoremstyle{remark}

\newcommand{\ft}{\widehat}
\newcommand{\mc}{\mathcal}

\newcommand{\Sp}{\mathbb{S}}

\newcommand{\C}{\mathbb{C}}

\newcommand{\R}{\mathbb{R}}

\newcommand{\bx}{\boldsymbol{x}}

\newcommand{\dx}{\text{\rm d}x}

\newcommand{\es}[1]{\begin{equation}\begin{split}#1\end{split}\end{equation}}

\newcommand{\ov}{\overline}
\renewcommand{\H}{\mc{H}}
\newcommand{\im}{{\rm Im}\,}

\newcommand{\bo}{\boldsymbol}

\newcommand{\vol}{{\rm vol}}

\begin{document}


\title[]{A Note on Band-Limited Minorants of an Euclidean Ball}
\author[]{Felipe Gon\c{c}alves}
\date{\today}
\subjclass[2010]{}
\address{University of Alberta, Mathematical and Statistical Sciences, CAB 632, Edmonton, Alberta, Canada T6G 2G1}
\email{felipe.goncalves@ualberta.ca}
\urladdr{sites.ualberta.ca/~goncalve}
\allowdisplaybreaks
\numberwithin{equation}{section}


\begin{abstract}
{We study} the Beurling-Selberg problem of finding band-limited $L^1$-functions that lie below the indicator function of an Euclidean ball. We compute the critical radius of the support of the Fourier transform for which such construction {can have a positive integral}.
\end{abstract}

\maketitle


\section{Introduction}
For a given $r>0$ we denote by $B^d(r)$ the closed Euclidean ball in $\R^d$ centered at the origin with radius $r>0$. We simply write $B^d$ when $r=1$. Define the following quantity
\begin{equation}\label{beta-def}
\beta(d,r)=\sup_{F} \int_{\R^d} F(\bx)dx,
\end{equation}
where the supremum is taken among functions $F\in L^1(\R^d)$ such that:
\begin{enumerate}
\item The Fourier Transform of $F(x)$, 
$$
\ft F(\xi) = \int_{\R^d} F(x)e^{2\pi ix\cdot \xi}\dx,
$$
is supported in $B^d(r)$;

\item $F(x)\leq {\bo 1}_{B^d}(x)$ for all $x\in \R^d$.
\end{enumerate}
We call such a function $\beta(d,r)$-admissible. A trivial observation is that $F\equiv 0$ is $\beta(d,r)$-admissible, hence $\beta(d,r)\geq 0$. Heuristically, such function $F(x)$ should exist and its mass should be close to $\vol(B^d)$ when $r$ is large. On the other hand, if $r$ is small, the mass of $F(x)$ should be close to zero and a critical $r_d>0$ should exist such that no function can beat the identically zero function for $r\leq r_d$. For this reason we define
\begin{equation*}
r_d=\inf \{r>0: \beta(d,r)>0\}
\end{equation*}
and it is this critical radius that we want to study in this paper.
\smallskip

The problem stated in \eqref{beta-def} has its origins with Beurling and Selberg which studied one-sided band-limited approximations for many different functions other than indicator functions with the purpose of using them to derive sharp estimates in analytic number theory (see the introduction of \cite{Va} for a nice first view). Although Selberg was one of the first to study the higher dimensional problem, it was first systematically analyzed by Holt and Vaaler in the remarkable paper \cite{HV}. They were able to construct non-zero {$\beta(d,r)$-admissible} functions for any $r>0$ and, most importantly, they established a fascinating connection of the $d$-dimensional problem with the theory of Hilbert spaces of entire functions contructed by de Branges (see \cite{dB}). They reduced the higher dimensional problem, after a radialization argument, to a weighted one-dimensional problem where the weight was given by a special function of Hermite-Biehler class, which in turn allowed them to use the machinery of homogeneous de Branges spaces to attack the problem. This new connection established by Holt and Vaaler started a new way of thinking about these kind of problems and ultimately inspired Littmann to completely solve the one-dimensional problem in \cite{L} by using a cleaver argument based on a special structure of certain de Branges spaces. Finally, using the ideas introduced by Littmann in \cite{L}, the problem of minorizing the indicator function of a symmetric interval was completely solved in \cite{CCLM} in the de Branges space setting.
\smallskip

{This paper was mainly motivated by the related problem where balls are substituted by boxes $Q(r)=[-r,r]^d$ and where practically nothing is known (see \cite{CEGK}). The box minorant problem is harder since it is a truly higher dimensional problem, whereas for the ball we can make radial reductions that transform it in a one-dimensional problem. Another interesting similar question, connected with upper bounds for sphere packings in $\R^d$, is studied in \cite{G} (see also \cite{CE}), where the author constructs a minorant $F(x)\leq {\bo 1}_{B^d}(x)$ with Fourier transform non-negative and supported in $B^d(j_{d/2,1})$ and such that it maximizes $\ft F(0)$ among all functions with these properties}\footnote{It is not the intention of this paper to give  a survey of related articles on the subject, which is very rich and full of subtleties, the purpose here is to {draw} a strait line between what he have so far and what we want to show.}.

\subsection{Main Result}

For any given parameter $\nu>-1$ let $J_\nu$ denote the classical Bessel function of the first kind. We also denote by $\{j_{\nu,n}\}_{n\geq 1}$ its positive zeros listed in increasing order. {The Bessel function of the first kind $J_\nu$ can be defined in a number of ways. We follow the treatise \cite{W} and define it for $\nu > -1$ and $\Re(z) >0$ by 
\begin{equation}\label{def_Bessel}
J_\nu(z) = \Big(\frac z 2\Big)^\nu \sum_{n=0}^{\infty} \frac{(-1)^n \big(\tfrac z 2\big)^{2n}}{n!\, \Gamma(\nu+n+1)}.
\end{equation}
For these values of $\nu$, one can check that on the half space $\{\Re(z) >0\}$ the Bessel functions defined by \eqref{def_Bessel} satisfy the differential equation
\begin{equation*}
z^2J''_\nu(z)+z J_\nu'(z)+(z^2-\nu^2)J_\nu(z)=0,
\end{equation*}
and that the following recursion relations hold
\begin{align*}
J_{\nu-1}(z)-J_{\nu+1}(z)&=2 J'_\nu(z),\\
J_{\nu-1}(z)+J_{\nu+1}(z)&=\frac{2\nu}{z} J_\nu(z).
\end{align*}
In particular we have $J_{-1/2}(z)=\sqrt{\frac{2}{\pi z}}\,\cos(z)$ and $J_{1/2}(z)=\sqrt{\frac{2}{\pi z}}\,\sin(z)$, which implies that
$j_{-1/2,1}=\pi/2$ and $j_{1/2,1}=\pi$.} The following is the main result of this paper. 

\begin{theorem}\label{crit-radius-thm}
We have 
$$
r_d = \frac{j_{d/2-1,1}}{\pi}.
$$
Moreover, if ${j_{d/2-1,1}} < \pi r < {j_{d/2,1}}$ then
$$
\beta(d,r)=\frac{(2/r)^d}{{|\Sp^{d-1}|}}\frac{\gamma_{\pi r}}{1+\gamma_{\pi r}/d},
$$
where $\gamma_{\pi r} = -\frac{\pi r J_{d/2-1}(\pi r)}{J_{d/2}(\pi r)}>0$. In particular we have
$$
\beta(d,r) = \frac{\pi^2 2^d}{r^{d-1}{|\Sp^{d-1}|}}\bigg(r-\frac{j_{d/2-1,1}}{\pi}\bigg) + O_d\bigg(r-\frac{j_{d/2-1,1}}{\pi}\bigg)^2
$$
for $r$ close to $\frac{j_{d/2-1,1}}{\pi}$.
\end{theorem}

\subsection*{{Remarks}} 
\begin{enumerate}
\item It is known that $j_{\nu,1} = \nu + 1.855757\nu^{1/3}+O(\nu^{-1/3})$ as $\nu\to\infty$ (see \cite[Section 1.3]{E}). This implies that $r_d = \frac{d}{2\pi} + \frac{1.855757}{2^{1/3}\pi}d^{1/3} + O(d^{-1/3})$ as $d\to\infty$. Heuristically, this means that if one wishes to non-trivially minorate (that is, beat the zero function) the indicator function of a ball of radius of order $\sqrt{d}$ then one needs frequencies of order at least $\sqrt{d}$. 

{\item The first $5$ values of $r_d$ rounded up to 4 significant digits are the following: $r_1=1/2$, $r_2 = 0.7655$, $r_3=1$, $r_4=1.220$ and $r_5=1.431$.}

\item Explicit expressions for $\beta(d,r)$ can also be tracked from \cite[Theorem 5]{CCLM}, but they involve sums of Bessel functions evaluated at Bessel zeros that can be quite complicated to grasp. Moreover, this is the case only when $\pi r$ is a zero of $J_{d/2-1}(z)$ or $J_{d/2}(z)$. If that is not the case, then writing a formula for $\beta(d,r)$ becomes pointless, since it will involve zeros of more complicated functions related to Bessel functions and this is not the purpose here.
\end{enumerate}

\section{Proof of Theorem \ref{crit-radius-thm}}

\subsection*{Step 1.} The first step is to reduce the higher dimensional by considering only radial functions. We can apply \cite[Lemmas 18 and 19]{HV} to reduce the $d$-dimensional problem to the following weighted one-dimensional problem
\es{\label{radial-id}
\beta(d,r) = \frac{|\Sp^{d-1}|}{2}\sup_{F} \int_{\R} F(x)|x|^{d-1}dx,
}
where $|\Sp^{d-1}|$ denotes the surface area of the unit sphere in $\R^d$ and the supremum is taken among functions $F\in L^1(\R,|x|^{d-1}\dx)$ such that:
\begin{enumerate}
\item $F(x)$ is the restriction to the real axis of an even entire function $F(z)$ of exponential type at most $2\pi r$, that is,
$$
|F(z)| \leq Ce^{2\pi r|\im z|}, \ \ \ z\in\C
$$
for some constant $C>0$;

\item $F(x)\leq {\bo 1}_{[-1,1]}(x)$ for all $x\in \R$.
\end{enumerate}

{In this framework the problem} becomes treatable with the theory de Branges spaces of entire functions. The latter generalize the well known Paley-Wiener spaces by using weighted norms given by  Hermite-Biehler functions. In what follows we briefly review the construction of a special family of de Branges spaces called {\it homogeneous spaces} which were introduced by de Branges (see \cite[Section 50]{dB} and \cite[Section 5]{HV}). We refer to \cite[Section 3]{HV} for a brief description of the general theory and also to \cite[Chapter 2]{dB} for the full theory. 

\subsection*{Step 2.} Let $\nu > -1$ be a parameter and consider the real entire functions\index{entire function} $A_\nu(z)$ and $B_\nu(z)$ given by
\begin{equation*}
A_{\nu}(z) = \sum_{n=0}^{\infty} \frac{(-1)^n \big(\tfrac12 z\big)^{2n}}{n!(\nu +1)(\nu +2)\ldots(\nu+n)} = \Gamma(\nu +1) \left(\tfrac12 z\right)^{-\nu} J_{\nu}(z)
\end{equation*}
and
\begin{equation*}
B_{\nu}(z) = \sum_{n=0}^{\infty} \frac{(-1)^n \big(\tfrac12 z\big)^{2n+1}}{n!(\nu +1)(\nu +2)\ldots(\nu+n+1)} = \Gamma(\nu +1) \left(\tfrac12 z\right)^{-\nu} J_{\nu+1}(z).
\end{equation*}
If we write
\begin{equation*}
E_\nu(z) = A_\nu(z) - iB_\nu(z),
\end{equation*}
then $E_\nu(z)$ is a Hermite--Biehler function, that is, it satisfies the following fundamental inequality
$$
|E_\nu(\ov z)| < |E_\nu(z)|
$$
for all $z\in\C$ with $\im z>0$. It is also known that this function does not have real zeros, that $E(iy)\in\R$ for all real $y$ (that is, $B_\nu(z)$ is odd and $A_\nu(z)$ is even), that $E_\nu(z)$ is of bounded type in the upper-half plane (that is, $\log |E_\nu(z)|$ has a positive harmonic majorant in the upper-half plane) and $E_\nu(z)$ is of exponential type $1$. We also have that
\begin{equation*}
 c|x|^{2\nu+1}\leq |E_{\nu}(x)|^{-2} \leq C |x|^{2\nu+1}
\end{equation*}
for all $|x|\geq 1$ and for some $c,C>0$. The homogeneous space $\H(E_\nu)$ is then defined as the space of entire functions $F(z)$ of exponential type at most $1$ and such that\footnote{ As an historical note, de Branges originally defined this space in another way but, in \cite[Lemma 16]{HV}, {the authors} showed that this is an equivalent definition.}
$$
\int_\R |F(x)|^{2}\,|E_{\nu}(x)|^{-2}\, {\\dx} < \infty.
$$
Using standard asymptotic expansions for Bessel functions one can show that $A_\nu,B_\nu\notin \H(E_\nu)$. As a particular case, observing that $E_{-1/2}=e^{-iz}$ we can deduce that $\H(E_{-1/2})$ coincides with the Paley-Wiener space of square integrable entire functions of exponential type at most $1$.

These spaces are relevant to our problem since we have the following magical identity
\begin{align}\label{magical-id}
a_\nu \int_\R |F(x)|^{2}\,|E_{\nu}(x)|^{-2}\, {\dx} =  \int_\R |F(x)|^2 \,|x|^{2\nu+1} \,{\dx}
\end{align}
for each $F\in\H(E_\nu)$, where $a_\nu = \,2^{2\nu+1}\, \Gamma(\nu+1)^{2}/\pi$. {For our purposes we will need an identity analogous of \eqref{magical-id}, but which allow us to compute integrals instead of $L^2$-norms. It can be derived as follows.} Let $F(z)$ be an entire function of exponential type at most $2$ such that $F(x)\leq {\bf 1}_{[-t,t]}(x)$ for some $t>0$ and $F\in L^1(\R,|x|^{2\nu+1}\dx)$. Since $G_n(x)=(\frac{\sin(x/n)}{x/n})^{n}$ belongs to $\H(E_\nu)$ for large $n$ and it converges to $1$ uniformly in compact sets as $n\to \infty$, we have that $4G_n(x)^2-F(x)\geq 0$ for all real $x$ (if $n$ is large and even) and this function has exponential type at most $2$. This implies that $4G_n(x)^2-F(x)=H_n(z)\ov{H_n(\ov z)}$ for all $z\in\C$, for some entire function $H_n(z)$ of exponential type at most $1$ (see \cite[Theorem 13]{dB}). We conclude that $H_n\in\H(E_\nu)$ and we obtain
\es{\label{magical-id-2}
a_\nu \int_\R F(x)\,|E_{\nu}(x)|^{-2}\dx & = a_\nu \int_\R \{4G_n(x)^2-|H_n(x)|^2\}\,|E_{\nu}(x)|^{-2}\dx \\ 
& = \int_\R \{4G_n(x)^2-|H_n(x)|^2\}|x|^{2\nu+1}\dx \\
& = \int_\R F(x)|x|^{2\nu+1}\dx.
}

\subsection*{Step 3.} Taking $\nu=d/2-1$, {we can apply the change of variables $x\mapsto x/(\pi r)$ in \eqref{radial-id} and use identity \eqref{magical-id-2}} to reduce the problem of minorizing the indicator function of an Euclidean ball to the following final form
$$
\beta(d,r) = \frac{(2/r)^d}{\pi{|\Sp^{d-1}|}}\Lambda_{E_{d/2-1}}^-(\pi r),
$$
where
$$
\Lambda_{E_{d/2-1}}^-(\pi r)=\sup_{F} \int_{\R} F(x)|E_{d/2-1}(x)|^{-2}\dx
$$
and the supremum is taken among functions $F\in L^1(\R,|E_{d/2-1}(x)|^{-2}\dx)$ such that:
\begin{enumerate}
\item $F(x)$ is the restriction to the real axis of an even entire function $F(z)$ of exponential type at most $2$;

\item $F(x)\leq {\bo 1}_{[-\pi r,\pi r]}(x)$ for all $x\in \R$.
\end{enumerate}
The above problem was completely solved in \cite{CCLM}. By all the previous discussion in Step $2$, we can apply \cite[Theorem 5 (i) and (iv)]{CCLM} to the function $E_{d/2-1}(z)$ (it actually can be applied to any $E_\nu(z)$) to derive that $\pi r_d=j_{d/2-1,1}$. Moreover, if $j_{d/2-1,1}<\pi r<j_{d/2,1}$ then \cite[Theorem 5 (iv)]{CCLM} also give us that
$$
\Lambda_{E_{d/2-1}}^-(\pi r) = \frac{\pi\gamma_{\pi r}}{1+\gamma_{\pi r}/d},
$$
where $\gamma_{\pi r} = -\frac{\pi r J_{d/2-1}(\pi r)}{J_{d/2}(\pi r)}>0$. A simple Taylor expansion leads to
$$
\Lambda_{E_{d/2-1}}^-(\pi r) = \pi^2 r(\pi r-j_{d/2-1,1})+O(\pi r-j_{d/2-1,1})^2
$$
and we finally obtain that
$$
\beta(d,r) = \frac{\pi^2 2^d}{r^{d-1}|{\Sp^{d-1}}|}\bigg(r-\frac{j_{d/2-1,1}}{\pi}\bigg) + O_d\bigg(r-\frac{j_{d/2-1,1}}{\pi}\bigg)^2.
$$

\end{document}